\documentclass[letterpaper, 10 pt, conference, onecolumn]{ieeeconf}  

\IEEEoverridecommandlockouts                              
\usepackage[loading]{tracefnt}

\usepackage[T1]{fontenc}
\usepackage{amsmath} 
\usepackage{amssymb}  
\usepackage{bbold,bbm}
\usepackage{siunitx}
\usepackage{algorithm, algorithmic, setspace}
\newcommand{\argmin}[1]{\underset{#1}{\operatorname{arg}\,\operatorname{min}}\;}

\newtheorem{remark}{Remark}
\newtheorem{example}{Example}
\newtheorem{assumption}{Assumption}

\newtheorem{theorem}{Theorem}
\newtheorem{proposition}{Proposition}

\newcommand{\R}{\mathbb{R}}
\newcommand{\B}{\mathcal{B}}
\renewcommand{\P}{\mathcal{P}}
\newcommand{\Q}{\mathcal{Q}}
\newcommand{\D}{\mathcal{D}}

\newcommand{\an}{\Tilde{n}} 
\newcommand{\am}{\Tilde{m}} 
\newcommand{\norm}[1]{\lVert#1\rVert}
\newcommand{\indic}{\mathbb{1}}

\title{\LARGE \bf Alternating direction method of multipliers\\ for polynomial optimization}

\author{V. Cerone \and S. M. Fosson \and S. Pirrera \and D. Regruto\thanks{$^{*}$ The authors are with the Dipartimento di Automatica e Informatica, Politecnico di Torino,
    corso Duca degli Abruzzi 24, 10129 Torino, Italy;
    e-mail: vito.cerone@polito.it, sophie.fosson@polito.it, simone.pirrera@polito.it, diego.regruto@polito.it;
}
}

\begin{document}
\maketitle
\thispagestyle{empty}
\pagestyle{empty}

\begin{abstract}
Multivariate polynomial optimization is a prevalent model for a number of engineering problems. From a mathematical viewpoint, polynomial optimization is challenging because it is non-convex. The Lasserre's theory, based on semidefinite relaxations, provides an effective tool to overcome this issue and to achieve the global optimum. However, this approach can be computationally complex for medium and large scale problems. For this motivation, in this work, we investigate a local minimization approach, based on the alternating direction method of multipliers, which is low-complex, straightforward to implement, and prone to decentralization. The core of the work is the development of the algorithm tailored to polynomial optimization, along with the proof of its convergence. Through a numerical example we show a practical implementation and test the effectiveness of the proposed algorithm with respect to state-of-the-art methodologies.
\end{abstract}

\section{Introduction}
Polynomial optimization problems (POPs) are concerned with the minimization of multivariate polynomial functions, over regions defined by polynomial equations and inequalities. Popular instances of POPs include linear, quadratic and mixed-integer programming, as well as partition problems on graphs, such as max-cut; see, e.g., \cite{lau09} for a list of examples. Beyond that, a wide class of engineering problems can be formulated via POPs. Among them, we mention set-membership estimation and identification \cite{cer12,cer14}, data-driven control \cite{cer17}, hybrid system identification \cite{fen10}, model predictive control \cite{can14,har16}, and optimal control \cite{rod22}. More recent applications include localization in narrowband internet-of-things \cite{sed21}, routing games \cite{yix22}, and estimation of Lipschitz constants in deep neural networks \cite{ton20}.

In general, POPs are non-convex and NP-hard, therefore the search of their global minima is challenging. To this purpose, solutions based on semidefinite programming (SDP) relaxations are proposed in the literature. In a nutshell, one can build a hierarchy of SDP problems that converge to the optimum under mild assumptions; see \cite{lasbook,par03,sho87} for a complete overview. The SDP approach recast POPs into convex optimization and this makes the problem affordable. The development of the related theory is mature and several software implementations are available, such as GloptiPoly \cite{hen02} and SparsePOP \cite{wak08}.

A drawback of the SDP approach is that the size of the involved matrices at $d$-degree of the hierarchy is proportional to $\left(\begin{array}{c}
                                                                                                n+d\\
                                                                                                    n
                                                                                                   \end{array}
\right),$   where $n$ is the number of variables of the original POP. This becomes computationally prohibitive for medium and large scale models, in particular when high $d$ are necessary to achieve the optimal solution.
To tackle this problem, in the recent literature,  SDP accelerating strategies are proposed, that exploit the presence of specific sparsity structures, such as the running intersection property \cite{wei18}, chordal sparsity \cite{fan19,zhe19}, or term sparsity \cite{wan21}. However, these sparsity structures are not always present in SDP relaxations of POPs, as discussed, e.g., in \cite{fan19}; therefore, these fast implementations are feasible only for some classes of POPs.
%


Given the described shortcomings of global optimization with SDP, one can relax POPs to the search of local minima or stationary points. In the literature, substantial attention is devoted to the local minimization of non-convex problems, e.g., via gradient-based methods, due also to the rise of deep learning techniques; we refer the reader to \cite{dan22} for a recent survey. Convergence of the algorithms is a delicate point in this framework, both in terms of guarantees and speed.

In this work, we propose a novel approach to local minimization of POPs, based on the alternating direction method of multipliers (ADMM, \cite{boyd11_admm}). ADMM is well-known in the context of convex optimization as an effective algorithm to minimize composite functionals, even in the presence of non-differentiable terms. ADMM structure is easy to implement and prone to decentralization. Its convergence is proven and analysed, see, e.g, \cite{boyd11_admm,hon17}. As to non-convex optimization, the convergence of ADMM is more critical to analyse and actually proven for some classes of problems; the main results are provided in \cite{hon16,wan19}.

The implementation of ADMM for POPs is not straightforward, since POPs are non-convex,  constrained problems, while ADMM is originally conceived for convex, unconstrained problems. The first contribution of this work is a suitable reformulation of POPs, on top of which we develop a novel ADMM-based procedure for POPs, that we name ADMM4POP. The second contribution is the proof of the convergence of ADMM4POP for the case of equality constraints, accounted for in a relaxed way. Furthermore, we illustrate the  implementation of ADMM4POP in some numerical examples and we compare its practical effectiveness to state-of-the-art global/local minimization methods.

The paper is organized as follows. In Sec. \ref{sec:problem}, we state the problem and we illustrate how to formulate it to apply ADMM. In Sec. \ref{sec:admm}, we illustrate the details of the proposed ADMM4POP algorithm, while in Sec. \ref{sec:proof} we prove its convergence. In Sec. \ref{sec:example}, we propose two numerical examples. Finally, in Sec. \ref{sec:concl} we draw some conclusions.
\section{Problem formulation}\label{sec:problem}
We consider a generic POP:
\begin{equation} \begin{aligned} \label{pop1}
    \min_{x\in\R^n} \quad & f(x) \\
    & \text{s.t.} \\
    & g_i(x) \leq 0 \quad i = 1,\dots,l \\
    & h_i(x) = 0 \quad i = 1,\dots,m
\end{aligned} \end{equation}
where $f,g_i,h_i$ are multivariate polynomials of the decision variable $x \in \R^n$; see, e.g., \cite{lasbook} for details.
The aim of this work is to develop an ADMM-based method for \eqref{pop1}. For this purpose, first of all, we need to rewrite \eqref{pop1} as an unconstrained, composite problem, with some specific properties. In this section, we illustrate the proposed reformulation.

The first step is to write \eqref{pop1} as a quadratic optimization problem (QOP). As discussed, e.g., in \cite{mev10,ell19}, any POP can be transformed into a QOP by adding suitable slack variables. Given a POP, different equivalent QOP representations are possible, as illustrated in \cite[Sec. 2.3]{mev10}.
In the rest of the paper, given two vectors $v,w$ of equal dimension, we write $v\leq w$ to denote the componentwise inequality. Moreover, we use the notation $A\succeq 0$ to indicate that a matrix $A$ is positive semidefinite.

For our purpose, we consider any transformation that provides a QOP of this form:
%
\begin{equation} \begin{aligned} \label{qop}
    \min_{x\in\R^{\an}}  \quad &  \dfrac{1}{2} x^\top A x + a^\top x\\
    & \text{s.t.} \\
    & Bx -b \leq 0 \\
    & Cx - c = 0 \\
    & x_ix_j = x_k, \quad (i,j,k) \in \B
\end{aligned} \end{equation}
where $x \in \R^{\an}$, $\an\geq n$, is the augmented vector of decision variables; $A \in  \R^{\an,\an}$ is symmetric and $A\succeq 0$;
$B \in \R^{l,\an}$, $C \in \R^{\am,\an}$, $a \in \R^{\an}$,  $b \in \R^l$, and  $c \in \R^{\am}$. We specify that $\tilde{m}$ is the number of equality constraints after the introduction of the slack variables.  Moreover, $\B\subseteq \{1,\dots,\an\}^3$ is the set of $3-$tuples $(i,j,k)$ that contains all the indices such that $x_ix_j = x_k$.


Before illustrating the transformation to obtain \eqref{qop}, we provide some remarks.
\begin{remark}\label{rem:A}
If the transformation returns a QOP where $A$ is not positive semidefinite, it is sufficient to introduce suitable quadratic constraints that reduce $ \dfrac{1}{2} x^\top A x + a^\top x$ to a linear function, i.e., $A = 0$.
\end{remark}
\begin{remark}\label{rem:B}
Without loss of generality, we assume that for each $(i,j,k)\in \B$, the indices $i$, $j$ and $k$ are mutually different. If not, it is sufficient to add new slack variables. For example, in the presence of the constraint $x_i^2 = x_k$, i.e., $(i,i,k)\in\B$, we can define a new variable $x_j$ and the constraint $x_j = x_i$ to obtain $x_ix_j = x_k$.
\end{remark}
\begin{remark}\label{rem:C}
Without loss of generality, we assume that each index $i,j,k \in \{1,\dots,\an\}$ appears at most once in each tuple $(i,j,k) \in \mathcal{B}$. This is not restrictive because if an overlap occurs, we can introduce an additional variable $x_h$, for each repeated variable $x_i$, together with constraints of the kind $x_h = x_i$.
\end{remark}

For simplicity, we describe the transformation algorithm from \eqref{pop1} to \eqref{qop} through an illustrative example.
\begin{example}
Let us consider the POP $$\min_{x\in\R^2}x_1^2 x_2.$$ By defining $x_3=x_1^2$, we obtain $\min_{x\in\R^2}x_3 x_2$ s.t. $x_1^2=x_3$.
Following Remark \ref{rem:A}, we notice that $A=\frac{1}{2}\left(\begin{array}{ccc}
                                                 0&0&0\\
                                                 0&0&1\\
                                                 0&1&0\\
                                                \end{array}\right)$ is not positive semi-definite. Therefore, we define $x_4=x_3x_2$.

Moreover, according to Remark \ref{rem:B},  we define $x_5=x_1$ to transform the constraint $x_1^2=x_3$ into $x_1x_5=x_3$.
In this way, $\B=\{(3,2,4),(1,5,3)\}$. To fulfill the non-overlap assumption described in Remark \ref{rem:C}, we add $x_6=x_3$.
In conclusion, we obtain the QOP
\begin{equation*}
 \begin{split}
  \min_{x\in\R^5}&\quad x_4\\
  &\quad\text{s.t. }\\
  &\quad x_5=x_1, \qquad x_6=x_3\\
  &\quad x_ix_j=x_k,~(i,j,k)\in\B=\{(3,2,4),(1,5,6)\}.
 \end{split}
\end{equation*}
\end{example}

\vskip0.8cm

Let us define the sets
\begin{equation} \label{setP}
    \P \doteq \{x\in\R^{\an} : x_i x_j = x_k, \quad \forall (i,j,k) \in \mathcal{B} \}.
\end{equation}
and
\begin{equation} \label{setD}
    \D \doteq \{x\in\R^{\an} | Bx -b \leq 0, Cx - c = 0\}.
\end{equation}
We notice that $\P$ is a non-convex semialgebraic set that describes all the non-convexity of the original problem. Instead, $\D$ is a convex set that includes all the linear constraints of \eqref{qop}. Given these definitions, we rewrite \eqref{qop} as
\begin{equation} \begin{aligned} \label{qop_admm_form}
    \min_{x\in\D }\quad &  \dfrac{1}{2} x^\top A x + a^\top x + \indic_{\P}(x)\\
\end{aligned} \end{equation}
where $\indic_{\P}(x)$ denotes the indicator function of $\P$, i.e., $\indic_{\P}(x)=0$ if $x\in\P$, and $\infty$ otherwise. 

\section{Proposed algorithm: ADMM4POP}\label{sec:admm}
In the previous section, we have shown that a generic  POP can be rewritten as a QOP of the form \eqref{qop_admm_form}. In this section, we prove that a POP formulated as in \eqref{qop_admm_form} is prone to be solved through ADMM.

In general, ADMM can be applied to composite problems of kind $\min_{x\in\R^n} F(x)+G(x)$; see \cite{boyd11_admm} for details and more general formulations. More specifically, to apply ADMM we split the problem as follows:
\begin{equation} \label{admm_form_generic} \begin{aligned}
    \min_{x,z\in\R^n} & \quad F(x) + G(z) \\
    & \quad\text{s.t.}   \qquad x=z
\end{aligned}\end{equation}
and we consider the augmented Lagrangian
\begin{equation}\label{auglag}
L_{\rho}(x,z)= F(x)+G(z)+\rho u^{\top}(x-z)+\frac{\rho}{2}\|x-z\|_2^2
\end{equation}
where  $\rho>0$ is a scalar penalty parameter and $u\in\R^{n}$ is the scaled dual variable; we refer the reader to \cite[Sec. 3.1.1.]{boyd11_admm} for details on the scaled ADMM.
In the following, we denote by $x_k$, $z_k$ and $u_k$ the estimates of the  variables $x$, $z$ and $u$, respectively, at a generic $k$-th iteration of the algorithm.

ADMM consists in a loop iterating over three steps:
\begin{subequations}
\begin{enumerate}
    \item $x$-update:
    \begin{equation} \label{x-update-generic}
    \begin{split}
        x_{\tau+1} &= \argmin{x\in\R^n}L_{\rho}(x,z)\\
        &=\argmin{x\in\R^n} F(x)+ \dfrac{\rho}{2} \norm{x-{z}_\tau+u_\tau}_2^2.
        \end{split}
    \end{equation}
    \item $z$-update:
    \begin{equation}  \label{z-update-generic}
     \begin{split}
        z_{\tau+1} &= \argmin{z\in\R^n}L_{\rho}(x,z)\\
        &=\argmin{z\in\R^n} G(z) + \dfrac{\rho}{2} \norm{x_{\tau+1}-{z}+u_\tau}_2^2.
        \end{split}
    \end{equation}
    \item $u$-update:
    \begin{equation}  \label{u-update-generic}
        u_{\tau+1} = u_\tau + x_{\tau+1}-z_{\tau+1}.
    \end{equation}
\end{enumerate}

\end{subequations}

Our aim is to apply ADMM to \eqref{qop_admm_form}. To this purpose,  we decompose \eqref{qop_admm_form}  as in \eqref{admm_form_generic}, by setting $F(x)= \dfrac{1}{2} x^\top A x + a^\top x$ and $G(z)=\indic_{\P}(z).$ Then, we have
\begin{equation} \label{admm_form_specific} \begin{aligned}
    \min_{x\in\D,z\in\R^{\an}} & \quad  \dfrac{1}{2} x^\top A x + a^\top x + \indic_{\P}(z) \\
    & \quad\text{s.t.}  \qquad x=z
\end{aligned}\end{equation}
and we implement the algorithm as in \eqref{x-update-generic}-\eqref{z-update-generic}-\eqref{u-update-generic}. While step \eqref{u-update-generic} is trivial, \eqref{x-update-generic} and \eqref{z-update-generic} require more computations.
\subsection{$x$-update step}\label{sec:xup}
According to \eqref{x-update-generic}, the $x$-update applied on \eqref{admm_form_specific} is
    \begin{equation} \label{xupdate}
        x_{\tau+1} = \argmin{x\in\D}  \dfrac{1}{2} x^\top A x + a^\top x + \dfrac{\rho}{2} \norm{x-z_\tau+u_\tau}_2^2.
    \end{equation}
This is equivalent to
    \begin{equation} \label{xupdate1} \begin{aligned}
        x_{\tau+1} &= \argmin{x\in\R^{\an}}  \dfrac{1}{2} x^\top \left(A+\rho I_{\an}\right) x + x^\top \left[a + \rho(u_\tau-z_\tau)\right] \\
        & \quad \text{s.t.} \\
        & \quad B x - b \leq 0, \\
        & \quad C x - c = 0
    \end{aligned}
    \end{equation}
    where $I_{\an}$ is the identity matrix of dimension $\an$.
    This is a convex quadratic program (QP) and can be solved efficiently to global optimality.

    Moreover, if no inequality constraints are present, the solution of \eqref{xupdate1} reduces to the solution of the linear system of the Karush-Kuhn-Tucker (KKT) conditions:
    \begin{equation}\label{KKT}
        \begin{bmatrix} A + \rho I & C^\top \\ C & 0 \end{bmatrix} \begin{bmatrix} x_{\tau+1} \\ \mu \end{bmatrix} = \begin{bmatrix} \rho(z_\tau - u_\tau) - a \\ c \end{bmatrix}
    \end{equation}
    where $\mu$ is the Lagrange multiplier; see \cite{bv04} for details. In sections \ref{sec:proof} and \ref{sec:example} we also discuss the relaxation of the equality constraints.
    %
%
\subsection{$z$-update step}\label{sec:zup}
According to \eqref{z-update-generic}, the $z$-update applied on \eqref{admm_form_specific} is
    \begin{equation} \label{zupdate} \begin{aligned}
        z_{\tau+1} &= \argmin{z\in\R^{\an}}\indic_{\P}(z) + \dfrac{\rho}{2}\norm{x_{\tau+1}-z+u_\tau}_2^2.
    \end{aligned}   \end{equation}
This is equivalent to
\begin{equation} \label{zupdate1} 
        z_{\tau+1} = \argmin{z\in\P} \norm{z- (x_{\tau+1}+u_\tau) }_2^2.
    \end{equation}
Since $\P$ is  non-convex, the solution to \eqref{zupdate1} is not straightforward. However, we can exploit the decoupling assumption discussed in Remark \ref{rem:C} to split the problem into smaller and affordable sub-problems. More precisely, on the one hand the quadratic cost in \eqref{zupdate1} is separable in its components; on the other hand, according to the decoupling assumption, we can partition $\P$ into separated subsets, each of them representing a constraint $x_i x_j = x_k$. Let $t=1,\dots,|\B|$ be the ordered indices of the elements of $\B$, i.e., each $(i,j,k)\in \B$ is labeled with a $t\in\{1,\dots,|\B|\}$. Then, we can write
\begin{equation}\label{partition}
 \begin{split}
  \P & = \P_1 \times \P_2 \times \dots \times \P_{|\B|}
 \end{split}
\end{equation}
where
\begin{equation}
\P_t  \doteq  \{ (x_i,x_j,x_k)  : x_ix_j = x_k \} \text{ for each } t=1,\dots,|\B|.
\end{equation}
Hence,
\begin{equation} \label{zupdate_split}
    z_{\tau+1}= \argmin{z\in\P_1 \times \dots  \times \P_{|\mathcal{B}|}} \sum_{t=1}^{|\mathcal{B}|} \norm{z_t - (x_{\tau+1}+u_\tau)_t }_2^2
\end{equation}
where $w_t \doteq (w_i,w_j,w_k)$ for any vector $w\in\R^{\an}$ and $(i,j,k)$  corresponding to the index $t$.

In conclusion, \eqref{zupdate_split} can be split into $|\mathcal{B}|$ sub-problems with only three variables, i.e.,
\begin{equation} \label{zupdate_splitted}
    (z_{\tau+1})_t = \argmin{z_t\in\P_t} \norm{z_t - (x_{\tau+1}+u_\tau)_t }_2^2.
\end{equation}
In other terms, $(z_{\tau+1})_t$ is the $\ell_2$ projection of  $v_t = (v_i,v_j,v_k) \doteq (x_{\tau+1}+u_\tau)_t \in \R^3$ onto $\P_t$. The problems in \eqref{zupdate_splitted} are non-convex; however, we can solve them in closed form. In fact, we have
\begin{equation} \label{zupdate_splitted_explicit} \begin{aligned}
    (z_{\tau+1})_t &= \arg\min_{z_t\in\R^3} (z_i - v_i)^2 + (z_j - v_j)^2 + (z_k - v_k)^2\\
    &\quad \text{s.t.} \qquad z_iz_j = z_k.
\end{aligned}   \end{equation}

Then, we plug $z_k= z_iz_j $ into the cost functional, and remove the  corresponding constraint:
\begin{equation} \label{zupdate_splitted_noconstr}
    (z_{\tau+1})_t = \arg\min_{z_t\in\R^3} (z_i - v_i)^2 + (z_j - v_j)^2 + (z_i z_j - v_k)^2.
\end{equation} %
Now,  we can explicitly find the minima by evaluating the first order conditions:  given $f_v(z_i,z_j) = (z_i - v_i)^2 + (z_j - v_j)^2 + (z_i z_j - v_k)^2$, we solve
\begin{equation} \label{gradient_zupdate}\begin{aligned}
    & \quad \nabla f_v(z_i,z_j) = 2\begin{pmatrix}
    z_i - v_i + z_i z_j^2 - z_j v_k \\
    z_j - v_j + z_i^2 z_j - z_i v_k
    \end{pmatrix} =  \begin{pmatrix} 0 \\ 0    \end{pmatrix}
\end{aligned}  \end{equation}
From the first equation, we obtain
\begin{equation} \label{z_i}
    z_i - v_i +z_iz_j^2 - z_jv_k =0 \Rightarrow z_i = \dfrac{v_i + z_jv_k}{1 + z_j^2 }\\
\end{equation}

By replacing \eqref{z_i} into the second equation of \eqref{gradient_zupdate}, we get
\begin{equation} \begin{aligned}
    & z_j - v_j + \left( \dfrac{v_i + z_jv_k}{1 + z_j^2} \right)^2z_j - \dfrac{v_i + z_jv_k}{1 + z_j^2}v_k = \\
    & = \dfrac{ z_j^5 + c_4z_j^4 + c_3z_j^3 + c_2z_j^2 + c_1z_j +c_0}{ \left(1 + z_j^2\right)^2 } = 0
\end{aligned}  \end{equation}
where
\begin{equation} \begin{aligned}
    & c_0 = - v_j - v_iv_k, \quad c_1 = v_i^2-v_k^2+1,  \\
    & c_2 = v_iv_k-2v_j, \quad c_3 = 2, \quad c_4 = - v_j.
\end{aligned} \end{equation}

As $1 + z_j^2\neq 0$, we obtain the candidate solutions by finding the zeros of a univariate polynomial of degree $5$, e.g., by computing the eigenvalues of the companion matrix
\begin{equation}
     \left[ \begin{array}{c|c}
    0 & I_{4}\\
    \hline
    -c_{0} & -c_{1} \ \dotsc \ -c_{4}
\end{array} \right]. \end{equation}
In this way, we get five candidate solutions, which are either real or complex and conjugate, and at least one of them is real. In turn, we evaluate the cost functional on each real candidate solution, which finally provides the desired minimum.

We notice that the $z$-update requires to compute the eigenvalues of a $5 \times 5$ matrix for each bilinear constraint, which may be burdensome if $|\B|$ is large. However, the problem can be fully parallelized, by exploiting the  separability of $\B$; therefore, in the presence of suitable hardware, the $z$-update can be solved very effectively.

\begin{algorithm}
\setstretch{1.2}
     \renewcommand{\algorithmicrequire}{\textbf{Input:}}
    \renewcommand{\algorithmicensure}{\textbf{Output:}}
  \caption{ADMM4POP}\label{alg:ADMM4POP}
  \begin{algorithmic}[1]

    %
    \STATE Initialization: $z_0,u_0\in\R^{\an}$
    \FORALL{$\tau=1,\dots,T_{stop}$}
    \STATE    $x_{\tau+1} = \argmin{x\in\D}  \dfrac{1}{2} x^\top A x + a^\top x + \dfrac{\rho}{2} \norm{x-z_\tau+u_\tau}_2^2,$
    via convex QP or KKT conditions, see Sec. \ref{sec:xup}
    %
    \STATE For each $t=1,\dots, |\B|$, $(z_{\tau+1})_t = \arg\min_{z_t\in\R^3} (z_i - v_i)^2 + (z_j - v_j)^2 + (z_i z_j - v_k)^2$
    by explicit computation, see Sec. \ref{sec:zup}
    \STATE $u_{\tau+1} = u_\tau + x_{\tau+1}-z_{\tau+1}$
    \ENDFOR
  \end{algorithmic}
\end{algorithm}
We summarize the proposed algorithm, denoted as ADMM4POP, in Algorithm \ref{alg:ADMM4POP}.

Concerning the stopping time $T_{stop}$, as suggested in \cite{boyd11_admm}, we stop the algorithm when the norms of both primal and dual residuals
$\norm{  x_\tau - z_\tau}$ and
     $\norm{\rho (z_{\tau-1} - z_{\tau})}$ are sufficiently small.

%

\section{Convergence analysis}
\label{sec:proof}
Proving the convergence of ADMM in non-convex problems is a challenging task. In \cite{li15}, the authors prove the convergence to a stationary point of single-block ADMM for a family of non-convex problems. Later, in \cite{hon16,wan19} extensions to multi-block problems are proposed. In \cite{hon16}, the terms of the cost functional that are  non-convex must be differentiable to prove the convergence. This is not the case of \eqref{qop_admm_form}, where  the indicator function of a non-convex set occurs, which is non-convex and non-smooth. Since in \cite{li15} the authors account for terms that are both non-convex and non-smooth, we leverage their results to develop the convergence analysis of ADMM4POP.

We remark that the analysis in \cite{li15} is developed for unconstrained problems; thus, for simplicity, in this work we recast problem \eqref{qop} into unconstrained optimization as well, while we leave the proof of convergence of the constrained version for future extended work. To this end, in \eqref{qop}, we assume that there are no inequality constraints $Bx\leq b$, while the equality constraints $Cx=c$ are accounted for in a relaxed way, i.e., by adding the term $\gamma\|Cx-c\|_2^2$ to the objective function, for some $\gamma>0$. Finally, we keep the bilinear constraints, because we represent them through the indicator function, see \eqref{qop_admm_form}.
In conclusion, we consider
\begin{equation} \begin{aligned} \label{relax_qop_admm_form}
    \min_{x\in\R^{\an}}  \quad & \dfrac{1}{2} x^\top A x + a^\top x+\gamma\|Cx-c\|_2^2+\indic_P(x).
\end{aligned} \end{equation}
As $\gamma\|Cx-c\|_2^2$ is quadratic and convex, \eqref{relax_qop_admm_form} has the form of \eqref{qop_admm_form}; then, we can apply ADMM4POP to it. We refer to this approach as relaxed ADMM4POP.

For completeness, we summarize a result from \cite{li15}, which is the basis for our proof of convergence.
%
%
Let us consider the composite problem \eqref{admm_form_generic} and its augmented Lagrangian \eqref{auglag}, under the following conditions.
\begin{assumption}\label{assu}$~~$

\begin{itemize}
\item[A1.] $F$ is semi-algebraic, twice continuously differentiable on $\R^n$ with a bounded, positive semidefinite Hessian $\nabla^2 F$;
\item[A2.] $G$ is a semi-algebraic, proper, lower semicontinuous function;
 \item[A3.] there exists $0<\zeta<\rho$ such that $F(x)-\frac{1}{2\zeta}\|\nabla F(x)\|^2$ is lower bounded;
 \item[A4.] the design parameter $\rho$ is designed so that $\rho I_n -\sqrt{2} \nabla^2F(x)\succeq 0$.
 \end{itemize}
\end{assumption}
We recall that a function is proper if it is finite somewhere and it does tend to $-\infty$; see, e.g., \cite[Sec. 2]{li15}.

\begin{theorem}\label{theo:1}\cite[Theorems 1-2-4]{li15}
Let us suppose that Assumption \ref{assu} holds.
Then, the sequence $(x_k,z_k,u_k)$ generated by ADMM converges to a point $(x^{\star},z^{\star},u^{\star})$, and $x^{\star}$ is a stationary point of $F(x)+G(x)$.
\end{theorem}

\begin{proposition}\label{cor:1}
Let us consider \eqref{relax_qop_admm_form}. The sequence $(x_k,z_k,u_k)$ generated by relaxed  ADMM4POP converges to a point $(x^{\star},z^{\star},u^{\star})$, and $x^{\star}$ is a stationary point of \eqref{relax_qop_admm_form}.
\end{proposition}

\proof
Let us consider $F(x)=\dfrac{1}{2} x^\top A x + a^\top x+\gamma\|Cx-c\|_2^2$ and $G(z)=\indic_P(z)$. In the following, we prove that all the points of Assumption \ref{assu} are fulfilled, which is sufficient to prove the thesis.

By construction of $A$, it is straightforward to prove that $\dfrac{1}{2} x^\top A x + a^\top x+\gamma\|Cx-c\|_2^2$ fulfills A1.

Regarding A2, by definition, $\indic_P$ is proper. Moreover, it is semi-algebraic, because $\P$ is a  semi-algebraic set, and the indicator function of a semi-algebraic set is semi-algebraic, see, e.g., \cite{att10}. Finally, the indicator function of a set is lower semi-continuous if the set is closed; then, we need to show that  $\P$ is closed. For this purpose, let $t$ be the index of $(i,j,k)$ as defined in Sec. \ref{sec:problem} and let us define
\begin{equation} \label{setPex}
    \Q_t \doteq \left\{ z \in \R^{\an} :  z_iz_j=z_k\right\}.
\end{equation}
We remark that $\Q_t$ is different from $\P_t$ defined in \eqref{partition}. In particular, it holds
    $\P= \bigcup_{t \in \mathcal{B}} \Q_t.$
%
Since
    $\Q_t = \left\{ z: -z_iz_j + z_k \leq 0 \right\} \cup  \left\{z : z_iz_j - z_k \leq 0 \right\}$
and since $z_iz_j - z_k$ and $-z_iz_j + z_k$ are continuous functions, their sub-level sets are closed. As a consequence $\Q_t$ is an intersection of closed sets, hence it is closed. In turn, $\P$ is a union of closed sets, hence it is closed. Then, we conclude that $\indic_{\P}$ is lower semi-continuous.

Finally, A3 holds because, in our problem, $F(x)-\frac{1}{2\zeta}\|\nabla F(x)\|^2$ is convex if we suitably design $\rho>\zeta>\|A+\gamma C^{\top}C\|_2$. In particular, if we set $\rho>\sqrt{2}\|A+\gamma C^{\top}C\|_2$, A4 is satisfied as well. \hfill$\blacksquare$

\section{Numerical experiments}
\label{sec:example}

In this section we propose two numerical experiments to test the effectiveness of the proposed ADMM4POP and to validate the theoretical convergence results. For ADMM4POP, we implement both the constrained version and the relaxed one, which is guaranteed to converge by Proposition  \ref{cor:1}. However, the constrained version as well is convergent in our experiments.

We compare ADMM4POP to two state-of-the-art methods, namely the interior-point method (IPM) and the SDP approach \cite{lasbook}.  IPM is a local algorithm like ADMM4POP. All the experiments are performed with MATLAB R2021b, on a PC with  AMD Ryzen 7 1700 8-Core CPU @ 3.00 GHz and 16 GB RAM.
We implement IPM through the function {\it{fmincon}} in the MATLAB Optimization Toolbox, while we use the SparsePOP package \cite{wak08} to implement the SDP method.

In the first experiment, we consider the following POP:
\begin{equation*} \begin{aligned}
    \min_{x \in \R^3} & \quad x_1^2x_2^2 + x_1^2 + q_1x_1 + x_2^2 + x_2x_3 +q_2x_2 + x_3^2 + q_3x_3 \\
    & \quad \text{s.t.}\\
    & \quad x_2x_3+x_1 = 10
\end{aligned} \end{equation*}
where $q_1\in[4,6]$, $q_2\in[-8,-6]$, and $q_3\in [1,3]$ are generated uniformly at random.
We reformulate the problem as in \eqref{qop}. First, we replace $x_4 = x_2x_3$ and $x_5=x_1x_2$.
Then, we define $x_6 = x_2$ to fulfill the splitting condition of Remark \ref{rem:C}. In conclusion, we obtain
\begin{equation*} \begin{aligned}
    \min_{x \in \R^6}\quad &  x_5^2 + x_1^2 + q_1 x_1 + x_2^2 + x_4 +q_2 x_2 + x_3^2 + q_3 x_3 \\
    &  \text{s.t.}\\
    &  x_4+x_1 = 10, \qquad  x_6=x_2\\
    &  x_ix_j = x_k,~ (i,j,k)\in \B = \{ (2,3,4), (1,6,5) \}
\end{aligned} \end{equation*}
which is in the form \eqref{qop_admm_form}. In particular, it is straightforward to verify that the quadratic cost functional is convex.

The stop threshold for ADMM4POP is set so that the obtained accuracy is comparable to that one of IPM. We randomly initialize all the variables with standard Gaussian distribution.  We set $\rho=2$, and $\gamma=10^3$ for relaxed ADMM4POP. In this experiment, SDP with relaxation order 2 always achieves the global minimum, thus it is used as reference.

In Table \ref{tab1}, we report the results over 500 random runs, in terms of accuracy and runtime statistics. We evaluate the error as the $\ell_2$ distance from the SDP solution, which is considered as ground truth. Both ADMM4POP and IPM substantially achieve the desired global minimum. As expected, in the relaxed version of ADMM4POP the error is slightly larger; however, the gap is not relevant in practice, that is, the global minimizer is identified with sufficient accuracy. As to the runtime, we can see that ADMM4POP is significantly faster than IPM and SDP.

\begin{table}[h]
    \centering
    \begin{tabular}{c|c| c c c}
    & Error & \multicolumn{3}{c}{Runtime ($\textrm{ms}$)} \\
    & mean & mean & min & max \\
    \hline
  SDP & -- &$63.5$&$59.0$&$88.5$\\
  IPM &$\SI{2.5e-05}{}$&$12.3$&$7.7$&$21.8$ \\
  ADMM4POP cns &$\SI{6.5e-05}{}$&$1.6$&$1.5$&$7.9$ \\
  ADMM4POP rel &$\SI{4.2e-04}{}$&$1.7$&$1.6$&$8.9$ \\
    \end{tabular}
    \caption{First experiment: constrained POP; statistics over 500 random runs.  ADMM4POP is implemented  in both constrained (cns) and relaxed (rel) versions.}
    \label{tab1}
\end{table}

In the second experiment, we deal with a system identification problem. We consider the  discrete-time dynamic system $y(k+1)=\alpha y(k)+\beta u(k+1)$, $k=0,\dots,K$, and we aim at recovering the parameters $\alpha$ and $\beta$ given the knowledge of the input $u\in\R^K$ and noisy measurements of the output $w\doteq y + \eta$, where $\eta\in\R^K$ represents the noise. In particular, we minimize the simulation error $\ell_2$-norm, i.e. we solve
\begin{equation*}
\begin{split}
\min_{y\in\R^K}\quad &\|y-w\|_2^2\\
&\text{s.t.}\\
&y(k+1)=\alpha y(k)+\beta u(k+1)~~k=1,\dots,K-1.
\end{split}
\end{equation*}
As illustrated in Sec. \ref{sec:problem}, we define the slack variables $\psi_k =\alpha$ and $\xi_k = \psi_k y(k)$ for $k=1,\dots,K-1$, to decouple the bilinear constraints. Thus, $\theta=(\alpha,\beta,y(1),\dots,y(K),\psi_1,\dots,\psi_{K-1},\xi_1,\dots,\xi_{K-1})\in\R^{2+K+2(K-1)}$ is the total vector of variables, and
$\B=\{(2+h,2+K+h,1+2K+h),h=1,\dots,K-1\}$.

\begin{table}[h!]
    \centering
    \centering
    \begin{tabular}{c|c| c c c}
    & Error & \multicolumn{3}{c}{Runtime ($\textrm{s}$)} \\
    & Mean & mean & min & max \\
    \hline
  SDP & $\SI{6.1e-04}{}$ &$10.16$&$5.87$&$22.99$\\
  IPM &$\SI{3.2e-04}{}$&$3.97$&$1.34$&$36.73$ \\
  ADMM4POP cns &$\SI{3.6e-04}{}$&$2.53$&$0.76$&$21.45$ \\
  ADMM4POP rel &$\SI{3.9e-04}{}$&$0.62$&$0.40$&$1.61$ \\
    \end{tabular}
    \caption{Second experiment: system identification; statistics over 500 random runs.  ADMM4POP is implemented  in both constrained (cns) and relaxed (rel) versions.}
    \label{tab:2}
\end{table}

We perform 500 random runs, with $\alpha,\beta\in (-1,-0.5)\cup(0.5,1)$ generated uniformly at random; $\eta$ is a white Gaussian noise with variance $10^{-4}$. The variables are initialized with Gaussian distribution. In each run, we measure $K=500$ output samples, the total number of variables being significantly larger with respect to the first experiment. For ADMM4POP, we set $\rho=1$ and $\gamma=10$. For the SDP approach, we set the relaxation order to 2, which is observed to be sufficient to achieve the global minimum.

In Table \ref{tab:2}, we collect the results. We define the error as the $\ell_2$ distance between the estimates and  the true parameters $\alpha,\beta$. All the considered algorithms identify the correct parameters, with similar accuracy. Regarding the runtime, ADMM4POP is faster than the other algorithms, with a more evident improvement in the relaxed version.

We remark that, in the proposed experiments, we perform the decoupled projections \eqref{zupdate_splitted} sequentially, which leaves room for further enhancement via parallelization, in particular in view of larger dimension problems.
\section{Conclusions}
\label{sec:concl}
In this work, we develop ADMM4POP, a low-complex strategy to implement the alternating direction method of multipliers for local minimization of polynomial optimization problems. By developing a suitable formulation of the problem, we prove that we can easily implement ADMM4POP for polynomial optimization, both in constrained and relaxed versions. Moreover,  we prove its convergence to a stationary point in the relaxed versions. Through  numerical examples, we show that ADMM4POP achieves the global minimum of polynomial problems with reduced run time with respect to state-of-the-art methods. As ADMM4POP is prone to decentralization, future work will envisage the parallelization of the proposed method and its application to large scale polynomial problems.


\bibliographystyle{IEEEtran}
\bibliography{optimization}

\end{document}